\documentclass[11pt,reqno]{amsart}

\usepackage{amssymb}
\usepackage[letterpaper,margin=1.18in]{geometry}
\usepackage[hidelinks]{hyperref}
\hypersetup{
  pdftitle={Rudin's Conjecture on Squares in Arithmetic Progressions from the abc Conjecture},
  pdfauthor={Alexey Kravatskiy}
}

\newtheorem{theorem}{Theorem}
\newcommand{\rad}{\operatorname{rad}}
\newcommand{\ii}{\mathrm{i}}
\newcommand{\Z}{\mathbb{Z}}

\title[Rudin's conjecture from $abc$]
 {Rudin's Conjecture on Squares in Arithmetic Progressions from the $abc$ Conjecture}
\author{Alexey Kravatskiy}
\address{MIRIAI (Moscow Independent Research Institute of Artificial Intelligence),
Moscow, Russia}
\email{kravatskii.a@miriai.org}
\date{November 27, 2026}

\begin{document}
\begin{abstract}
Assuming the classical $abc$ conjecture, we prove that every arithmetic
progression of $N$ integers with nonzero common difference contains
$O(\sqrt N)$ squares, uniformly in its initial term and common difference.
From four representations of the same integer as a sum of two squares,
Gaussian factorization and the Pl\"ucker relation give an equation involving
the fourth powers of three pairwise coprime integers, with coefficients
bounded polynomially in $N$. When a primitive progression contains more
than $5\sqrt N$ positive squares, a weighted cut argument combines the
equations obtained from five representations to bound the common difference
by a power of $N$. Finally, considering the square roots of the initial term
modulo the common difference gives the uniform square-root bound.
\end{abstract}
\maketitle

\section{Introduction}

For integers $N\ge1$, $a$, and $q\ne0$, let $Q(N;q,a)$ be the number
of squares among $a,a+q,\ldots,a+(N-1)q$, and set
$Q(N)=\max_{a,q:\,q\ne0}Q(N;q,a)$.
Rudin's conjecture asserts that $Q(N)\ll\sqrt N$ \cite{Rudin,CG}.
This order is sharp, as the progression $1,2,\ldots,N$ shows.

For a positive integer $n$, put $\rad(n)=\prod_{p\mid n}p$.
We use the usual $abc$ inequality in the following form: for a fixed
$\varepsilon>0$ there is a constant $C_\varepsilon$ such that
\begin{equation}\label{eq:abc}
 W\le C_\varepsilon\,\rad(UVW)^{1+\varepsilon}
\end{equation}
for all coprime positive integers $U,V,W$ with $U+V=W$.

\begin{theorem}\label{thm:main}
If \eqref{eq:abc} holds for a single fixed $0<\varepsilon<1/3$, then
$Q(N)\ll_\varepsilon\sqrt N$.
\end{theorem}

In particular, no $N^\varepsilon$ loss occurs in the exponent of the
conclusion. The theorem does not settle the stronger conjecture, studied by
Gonz\'alez-Jim\'enez and Xarles \cite{GJX}, that $Q(N)=Q(N;24,1)$ for
$N\ge6$.

\medskip
\noindent\textbf{Lean verification.}
The argument at $\varepsilon=1/4$ has been verified in Lean~4
\cite{RudinLean}. The formal theorem takes $abc$ as an explicit hypothesis.

\medskip
\noindent\textbf{Related work.}
Bombieri, Granville and Pintz \cite{BGP} obtained the unconditional bound
$Q(N)\ll N^{2/3}(\log N)^{O(1)}$; Bombieri and Zannier \cite{BZ}
reduced the power of $N$ to $3/5$. Bourgain and Demeter \cite{BD} obtained
the square-root bound when the common difference has a bounded number of
divisors. For polynomially bounded differences, Novakovi\'c
\cite[Theorem~1.4]{Novakovic} obtained the exponent $1/2+\eta$ for every
fixed $\eta>0$. Theorem~\ref{thm:main} assumes \eqref{eq:abc} but imposes
no restriction on the common difference and has no logarithmic loss.

\medskip
\noindent\textbf{Relation to earlier approaches.}
Cilleruelo and Granville \cite[Section~5, Theorem~12]{CG} reduce the problem
to representations as sums of two squares. They use Gaussian divisibility
and a conjecture on lattice points on circles to bound the common difference,
then use a quadratic congruence to bound the number of squares. In a
separate argument \cite[Section~9]{CG}, they use a strengthening of $abc$
due to Baker on a relation whose degree depends on the gaps between five
squares, obtaining a bound exponential in a power of $N$.
Croot, Mao and Yip
\cite[Theorem~2.5 and Remark~2.6]{CMY} obtained almost sharp conditional
bounds for powers $k\ge4$ and did not treat $k=2,3$ in that theorem.

Our main step is to derive a polynomial bound on the common difference
from $abc$ within the circle reduction. We apply the Pl\"ucker relation
to each set of four out of five representations of the same integer as a
sum of two squares. Gaussian factorization and rounding of valuations
give an equation with three fourth-power terms. Their bases are pairwise
coprime and their coefficients are bounded polynomially in $N$. Up to
polynomial factors in $N$, the largest base contributes its fourth power
to the size of the largest term, whereas the three bases contribute at
most its third power to the radical. The $abc$ inequality therefore
bounds each base. Each partition of the five representations into parts
of sizes two and three restricts to a balanced partition on at least one
of the chosen sets of four. The bounds for the bases, together with
estimates for the norms of the Gaussian gcds, give a polynomial bound on
the circle's squared radius and hence on the common difference. The
quadratic congruence then gives the square-root estimate.

\section{Five representations on a circle}

Reverse the progression if necessary to make $q>0$. We may assume it
contains a square. Write $\gcd(a,q)=b^2g$ with $g$ squarefree and put
$a=b^2g A$, $q=b^2g Q$, where $\gcd(A,Q)=1$.
If $a+qn=x^2$, then $bg$ divides $x$ and $g$ divides $A+Qn$.
Existence of such an index implies $\gcd(g,Q)=1$. The admissible indices
thus lie in one residue class modulo $g$. Let $n_0$ be the smallest
index in this class in the original progression. Dividing the terms
with indices $n=n_0+gk$ by $(bg)^2$ gives the progression
\[
 \frac{A+Qn_0}{g}+Qk.
\]
It has length at most $N$ and contains the same number of squares.
Its initial term is coprime to $Q$ because $A+Qn_0$ is. We may therefore assume
$\gcd(a,q)=1$. The term zero contributes at most one square and is
ignored below.

Suppose there are $m>5\sqrt N$ positive squares. The $\binom m2$ unordered
pairs of distinct positive square terms have at most $2N-1$ possible
index sums. Since $\binom m2\ge m^2/3>4(2N-1)$, some index sum $t$
occurs for at least five pairs.
Write the corresponding representations as
\begin{equation}\label{eq:circle}
 \begin{gathered}
 S=x_i^2+y_i^2=2a+qt,\quad x_i>y_i>0,
 \\
 d_i=\frac{x_i^2-y_i^2}{q}\in\{1,\ldots,N-1\},
 \qquad 1\le i\le5.
 \end{gathered}
\end{equation}
The $d_i$ are distinct and $d_i\equiv t\pmod 2$.
Order $z_i=x_i+\ii y_i$ by increasing argument, so that $d_i>d_j$
when $i<j$.

\section{Gaussian factorizations of the minors}

Put $a_i=qd_i$ and $b_i=2x_i y_i$, so that $z_i^2=a_i+\ii b_i$ and
$a_i^2+b_i^2=S^2$. For $i<j$ write
\[
 D_{ij}=a_i b_j-b_i a_j,\qquad E_{ij}=a_i b_j+b_i a_j.
\]
Both are positive. Let $g_{ij}$ be a Gaussian gcd of $z_i,z_j$, and put
$u=z_i/g_{ij}$, $v=z_j/g_{ij}$. They are coprime and have equal norm.
Since $u\mid v\bar v$ and $\gcd(u,v)=1$ in $\Z[\ii]$, we have
$u\mid\bar v$. Equality of norms then shows that $v$ is a unit
multiple of $\bar u$.
Set
\[
 K'_{ij}=|g_{ij}|^2,\qquad K_{ij}=|u|^2,\qquad K_{ij}K'_{ij}=S .
\]
The product $\bar z_i z_j=K'_{ij}\bar u v$ is divisible by $K'_{ij}$
in $\Z[\ii]$. Also $z_i z_j=g_{ij}^2uv$ is divisible by $K_{ij}$,
since $uv$ is a unit multiple of $u\bar u$. Squaring and taking
imaginary parts gives $(K'_{ij})^2\mid D_{ij}$ and $K_{ij}^2\mid E_{ij}$.

In fact
\begin{equation}\label{eq:J}
 J^-_{ij}=\frac{D_{ij}}{q(K'_{ij})^2},\qquad
 J^+_{ij}=\frac{E_{ij}}{qK_{ij}^2}
\end{equation}
are integers. Indeed, $q$ divides both minors and
$\gcd(q,S)=\gcd(q,2a)=\gcd(q,2)$.
If $q$ is even then $x_i,y_i$ are all odd, so $v_2(S)=1$,
$v_2(K'_{ij})=1$, and $v_2(K_{ij})=0$. Now
$D_{ij}=2q(d_i x_j y_j-d_j x_i y_i)$ and the expression in parentheses
is even since $d_i\equiv d_j\pmod2$. This supplies the additional
factor of four required for $J^-_{ij}$.

The identity $D_{ij}E_{ij}=q^2 S^2(d_i^2-d_j^2)$ gives
\begin{equation}\label{eq:J-bound}
 J^-_{ij}J^+_{ij}=d_i^2-d_j^2,\qquad
 1\le J^\pm_{ij}\le N^2.
\end{equation}
Also $D_{ij}\le a_i b_j\le qNS$. By \eqref{eq:J} we obtain
\begin{equation}\label{eq:K-bound}
 (K'_{ij})^2\le NS,\qquad K_{ij}\ge\sqrt{S/N}.
\end{equation}

\section{A quartic equation from four representations}

Choose any four of the five indices and relabel them in angular order.
The Pl\"ucker relation for the minors of the vectors $(a_i,b_i)$ reads
\[
 D_{12}D_{34}+D_{14}D_{23}=D_{13}D_{24}.
\]
Number the three pairings in the order of these terms. For a pairing
$\nu=(ij,kl)$ define $X_\nu=K'_{ij}K'_{kl}$ and $A_\nu=J^-_{ij}J^-_{kl}$.
After cancelling $q^2$ the relation becomes
\begin{equation}\label{eq:plucker}
 A_1X_1^2+A_2X_2^2=A_3X_3^2,\qquad 1\le A_\nu\le N^4.
\end{equation}

For each odd prime $p$ that splits in $\Z[\ii]$ and satisfies
$p^e\Vert S$, fix a Gaussian prime $\pi$ above $p$.
Let $\alpha_i=v_\pi(z_i)$, so that
$v_{\bar\pi}(z_i)=e-\alpha_i$. Choose $\beta_i\in\{0,1\}$ with
$e\beta_i$ nearest to $\alpha_i$, choosing a tie once for each $p,i$
throughout the proof. Then
\begin{equation}\label{eq:valuation}
 |\alpha_i-e\beta_i|\le\tfrac12 v_p(d_i),\qquad
 v_p(K'_{ij})=e-|\alpha_i-\alpha_j|.
\end{equation}
For the first inequality, $p^{\min(\alpha_i,e-\alpha_i)}$ divides both
$x_i$ and $y_i$, so its square divides $qd_i$; here $p\nmid q$.
For the second, take the smaller
valuation at each of $\pi,\bar\pi$ in the gcd.

At the prime $p$, the rounded valuation of $X_\nu$ obtained by
replacing $\alpha_i$ with $e\beta_i$ is
$e(\mathbf{1}_{\beta_i=\beta_j}+\mathbf{1}_{\beta_k=\beta_l})$.
The error is at most $\tfrac12\sum_i v_p(d_i)$ by \eqref{eq:valuation}.
For a $2+2$ split of the labels, the rounded valuations of the
three $X_\nu$ are $2e,0,0$ in some order; the valuation $2e$
belongs to the pairing that matches equal labels. In all other cases
the three rounded valuations agree. Let $P_\nu$ be the product of
the $p^e$ whose labels form a $2+2$ split with equal labels paired by $\nu$.
The three integers $P_\nu$ are pairwise coprime. At inert primes and at
$2$ the actual valuations of all three $X_\nu$ agree: a Gaussian prime
fixed up to units by conjugation has the same valuation in every $z_i$
of norm $S$.

Let $H=\gcd(X_1,X_2,X_3)$ and $D_0=\prod_i d_i\le N^4$. Subtracting
the minimum valuation from each of the three valuations changes the
rounding error by at most another $\tfrac12\sum_i v_p(d_i)$.
Consequently
\[
 \frac{X_\nu}{H}=P_\nu^2 t_\nu
\]
where $t_\nu>0$ is rational and its reduced numerator and denominator
each divide $D_0$. Divide \eqref{eq:plucker} by $H^2$ and multiply by
$D_0^2$. Since $D_0 t_\nu$ is an integer at most $D_0^2$, we obtain
\begin{equation}\label{eq:quartic}
 c_1P_1^4+c_2P_2^4=c_3P_3^4,\qquad
 c_\nu=A_\nu(D_0t_\nu)^2\in\Z,\qquad 1\le c_\nu\le N^{20}.
\end{equation}

Put $M=\max_\nu P_\nu$. The gcd of the three terms in
\eqref{eq:quartic} divides $c_1c_2c_3$, since a prime divides at most
one $P_\nu$. After division by the gcd the terms are pairwise coprime,
since the first two sum to the third. The largest term is at least
$M^4/N^{60}$, and the radical of the product is at most $N^{60}M^3$.
Applying \eqref{eq:abc} gives
\[
 M^{1-3\varepsilon}\le C_\varepsilon N^{120+60\varepsilon}.
\]
In particular, for each of the five possible quadruples,
\begin{equation}\label{eq:P-bound}
 \log P_\nu\ll_\varepsilon\log(2N).
\end{equation}

\section{Weighted cuts}

Restore all five indices. For each odd split prime $p$ with $p^e\Vert S$,
the labels $\beta_i$ define a cut of weight $e\log p$. For a pair
$i,j$, the total weight of the cuts separating $i$ and $j$ differs
from $\log K_{ij}$ by at most $(\log d_i+\log d_j)/2$, by
\eqref{eq:valuation}. Inert primes and $2$ do not contribute to
$K_{ij}$. Summing over the ten pairs and using \eqref{eq:K-bound}
shows that the sum of the ten separation weights is at least
$5\log S-15\log N$.

Let $T$ be the total weight of the cuts with parts of sizes two and
three. Such a cut separates six pairs; any other cut separates at
most four. Since the total weight is at most $\log S$,
\begin{equation}\label{eq:cuts}
 5\log S-15\log N\le 4\log S+2T.
\end{equation}
Every cut with parts of sizes two and three restricts to a $2+2$ cut
on at least one of the five quadruples. Thus $T$ is at most the sum
of the logarithms in \eqref{eq:P-bound} over all five quadruples
and all three pairings in each. It follows from \eqref{eq:cuts} that
$\log S\ll_\varepsilon\log(2N)$. Since $qd_i<S$, there is a fixed
$B=B(\varepsilon)$ such that $q\ll_\varepsilon N^B$.

\section{The quadratic congruence}

The positive square roots of the square terms lie in an interval of length
at most $\sqrt{qN}$, since the difference between the largest and smallest
square terms is less than $qN$. They satisfy $x^2\equiv a\pmod q$.
Because $\gcd(a,q)=1$, there are at most two solutions modulo each odd
prime power dividing $q$ and at most four modulo the power of $2$.
The Chinese remainder theorem therefore gives
\[
 m\ll 2^{\omega(q)}\left(1+\sqrt{N/q}\right),
\]
where $\omega(q)$ is the number of distinct prime divisors of $q$.
If $q\le N$, the bound $2^{\omega(q)}\ll\sqrt q$ gives
$m\ll\sqrt N$. If $N<q\ll_\varepsilon N^B$, fix
$0<\eta<1/(2B)$ and apply $2^{\omega(q)}\ll_\eta q^\eta$.
Both divisor estimates follow by separating the finitely many primes
with $p^\eta<2$ from the rest. This proves $m\ll_\varepsilon\sqrt N$
also in the second range. Together with the case $m\le5\sqrt N$ and
the initial normalization, this proves Theorem~\ref{thm:main}.
\qed

\section*{Acknowledgments}
The author posed the problem and asked whether the Riemann hypothesis
could be used to prove Rudin's conjecture. OpenAI's Dot, using GPT-6.2
Astra with Ultra reasoning effort, proposed instead to derive Rudin's
conjecture from the classical \(abc\) conjecture and developed the
conditional proof presented here. At the author's request, Dot also
formalized the same argument in Lean 4.

\bibliographystyle{amsplain}
\bibliography{references}
\end{document}